\catcode`\@=11

\font\tenmsa=msam10 \font\sevenmsa=msam7 \font\fivemsa=msam5
\font\tenmsb=msbm10
\font\sevenmsb=msbm7 \font\fivemsb=msbm5 \newfam\msafam
\newfam\msbfam
\textfont\msafam=\tenmsa \scriptfont\msafam=\sevenmsa
\scriptscriptfont\msafam=\fivemsa \textfont\msbfam=\tenmsb
\scriptfont\msbfam=\sevenmsb \scriptscriptfont\msbfam=\fivemsb

\def\hexnumber@#1{\ifnum#1<10 \number#1\else \ifnum#1=10
A\else\ifnum#1=11
 B\else\ifnum#1=12 C\else \ifnum#1=13 D\else\ifnum#1=14
E\else\ifnum#1=15
 F\fi\fi\fi\fi\fi\fi\fi}

\def\msa@{\hexnumber@\msafam} \def\msb@{\hexnumber@\msbfam} 

\mathchardef\hbar="0\msb@7E
\mathchardef\blacktriangleright="3\msa@49
\mathchardef\blacktriangleleft="3\msa@4A

\def\Bbb{\ifmmode\let\next\Bbb@\else
\def\next{\errmessage{Use \string\Bbb\space only in math
mode}}\fi\next}
\def\Bbb@#1{{\Bbb@@{#1}}} \def\Bbb@@#1{\fam\msbfam#1}

\catcode`\@=\active

\def\R{{\Bbb R}}
\def\C{{\Bbb C}}

\def\CR{\hbox{{$\cal R$}}}  
\def\CH{\hbox{{$\cal H$}}}
\def\vecu{{\bf u}}
\def\vecx{{\bf x}} 

\def\grav{{\scriptscriptstyle G}}

\def\del{\partial}

\def\eps{{\epsilon}}
\def\lcross{{>\!\!\!\triangleleft}}
\def\cobicross{{\triangleright\!\!\!\blacktriangleleft}}
\def\bicross{{\blacktriangleright\!\!\!\triangleleft}}
\def\tens{\mathop{\otimes}}
\def\isom{{\cong}}
\def\Ad{{\rm Ad}}

\def\<{\langle}
\def\>{\rangle}
\def\o{{}_{\scriptscriptstyle(1)}}
\def\t{{}_{\scriptscriptstyle(2)}}

\def\eqn#1#2{\begin{equation}#2\label{#1}\end{equation}}
\def\ceqn#1#2{\begin{equation}\label{#1}
\begin{array}{c}#2\end{array}\end{equation}}

\documentstyle[11pt]{article}
\newtheorem{lemma}{Lemma}[section]
\newtheorem{propos}[lemma]{Proposition}

\begin{document}
\begin{center} {\LARGE QUANTUM GEOMETRY AND THE PLANCK SCALE}\\  
 \baselineskip 13pt{\ }
{\ }\\
S. Majid\footnote{Royal Society University Research Fellow and Fellow
of
Pembroke College, Cambridge}\\
{\ }\\
Department of Mathematics, Harvard University\\
Science Center, Cambridge MA 02138, USA\footnote{During the calendar
years 1995
+ 1996}\\
+\\
Department of Applied Mathematics \& Theoretical Physics\\
University of Cambridge, Cambridge CB3 9EW\\
\end{center}

\begin{center}
December 1996
\end{center}
\vspace{10pt}
\begin{quote}\baselineskip 13pt 
We consider some general aspects of the new noncommutative or
quantum geometry coming out of the theory of quantum groups, in
connection with  Planck scale physics. A generalisation of Fourier
or wave-particle duality on curved spaces emerges. Another feature
is the need for particles with fractional or braid statistics. The
conformal group also has a special role.  
\end{quote}
\baselineskip 21pt
\section{Introduction}

In this contribution I would like to discuss some of the conceptual
and physical issues surrounding the approach to noncommutative
geometry coming out of quantum groups and braided groups, in
somewhat greater depth than I had time for in my lecture in Goslar.
The discussion is intended to be intelligible to non-specialists and
may hopefully serve as an invitation to the field. Technical
material including recent results in quantum and braided geometry may
be found in my contribution to the companion `Quantum Groups' volume
of these Proceedings. Basic material can also be found in my textbook
on quantum groups\cite{Ma:book}.

The need for some kind of quantum geometry has
been clear enough since the birth of quantum mechanics itself: how to
extend ideas of gauge theory, curvature and non-Euclidean geometry to
the the situation where coordinates are noncommuting operators. If
one ever wants to unify quantum mechanics and gravity into a single
exact theory then this is probably a prerequisite.  We begin, in
Section~2, by postulating some fundamental features which any
satisfactory such quantum geometry should have. One of them is an
extension of wave-particle duality to curved
space\cite{Ma:pri}\cite{Ma:pla}. 

In Section~3 we consider quantum groups as examples of quantum
geometry. Actually, this is not at all the context in which the more
famous $q$-deformed enveloping algebras $U_q(g)$ arose (which is that
of `generalised symmetries' of exactly solvable lattice models; there
$q$ is an anisotropy parameter and not directly related to Planck's
constant). However,  at about the same time as these $U_q(g)$ were
being introduced in the mid 1980's, another completely different
class of quantum groups $\C(G_1)\bicross U(g_0)$ was being introduced
(by the author) in another context, which was precisely the context
of
Planck scale physics and an algebraic approach to quantum-gravity.
These {\em bicrossproduct} quantum groups really do arise as quantum
algebras of observables\cite{Ma:pla}. Here $G_0G_1$ is a
Lie group factorisation and $g_0$ is the Lie algebra of $G_0$. 
More recently,  gauge theory etc. has been developed for quantum
geometry\cite{BrzMa:gau}\cite{Ma:diag}.  

As soon as spaces become noncommutative or `quantum', their
symmetries naturally become generalised as quantum groups too. Here
the quantum groups $U_q(g)$ play a more central role. It turns out as
a new feature of quantum geometry that anything on which a quantum
group acts acquires braid statistics. In other words. not only are
algebras noncommutative but tensor products become noncommutative
too.
This means that at the Planck scale one should expect not only bosons
and fermions but complicated braid statistics as well. This gives a
systematic `braided approach'\cite{Ma:introp} to quantizing
everything
uniformly, and is the topic of Section~4. 

Obviously everything we may want to say about the Planck scale here
will be speculative. However, mathematics can tell us that certain
assumptions will force us to certain conclusions on the mathematical
structure, without yet knowing realistic models.
Moreover, quantum geometry is probably necessary not only at
the Planck scale or in quantum cosmology, but also in the resolution
of those
paradoxes in quantum mechanics which are characterised by a conflict
between the macroscopic geometry of measuring apparati and quantum
mechanical evolution; it should provide the right language to
correctly formulate such questions.

Note also that most discussions of Planck scale physics, including
string theory, work within an underlying classical geometry (e.g.
inside a path integral.) This is practical
but not really justified: classical geometry should emerge from
a deeper quantum world and not vice-versa. Even Professor Fredenhagen
in his talk {\em assumed exact classical Poincar\'e symmetry} without
any real justification, except that this is a necessary assumption to
be able apply the known methods of algebraic quantum field theory.
This is still `looking for the key under the streetlight'; by
contrast the quantum geometry 
programme advocated here seeks primarily to understand first what
mathematically natural quantum geometry is out there, before making
predictions. The success of General Relativity can be attributed in
part to the fact that Einstein already had a fairly complete
mathematical theory of Riemannian geometry to bring to bear.  
For example, modified uncertainty relations based on modified
commutation relations $[x,p]$ are only as meaningful as the
justification that a particular operator $x$ be called position and 
a particular operator $p$ be called momentum, which  should generally
come from a quantum geometrical picture.

\section{Inventing quantum geometry}

Riemannian geometry as we know it arose in two stages. First, one can
work extrinsically with surfaces embedded in Euclidean space, or
submanifolds $M\subset \R^n$. This is a Gaussian approach. Riemann 
realised, however, that one needs also a more intrinsic notion of
geometry determined by structure on $M$ itself, particularly
since actual space or (after Einstein) spacetime is what we directly
experience rather than any embedding space. This leads to our modern
form of non-Euclidean or Riemannian geometry. The situation with
quantum mechanics can be viewed analogously. The role of $\R^n$ is
played in a certain sense by $B(\CH)$, the operators on a Hilbert
space. However, the real observables in the system are usually some
$*$-subalgebra $A\subset B(\CH)$. We broaden the
problem a little and consider the intrinsic structure of quite
general algebras $A$ (not only normed $*$-algebras over $\C$).  In my
opinion
we are in a similar situation to Riemann: how to develop a language
to describe the intrinsic  geometric structure on an algebra. Note
that for any manifold we can consider the algebra of smooth functions
on it and formulate our usual geometrical notions in those terms; the
key difference in quantum geometry is that we do not want to be
limited to commutative algebras.

Also, let us say from the start that the real physical motivation
for quantum geometry applies directly only to phase
spaces, which lose their points on quantisation due to the
uncertainty principle; the coordinates of phase space no longer fully
commute. However, a generalised framework for phases spaces probably
entails developing a generalised framework for manifolds in general.
Apart from this zeroth assumption, we postulate also:

\begin{enumerate}
\item {\bf Richness} -- at least as much `flabbiness' in the variety
of examples and structures (gauge fields, curvature etc.) as
classically.
\item {\bf Quantization without classical assumptions} -- classical
geometry should emerge as a possible limit, not be built in from the
start by assuming a Poisson manifold.
\item {\bf Uniformity of quantisation} -- the whole geometrical
`zoo' should be quantised together, coherently, not one object at a
time.
\item {\bf Duality between geometry and quantum mechanics} -- wave
particle duality should be maintained in some form.
\end{enumerate}

Axiom~1 here is not as empty as it may seem. There are
many noncommutative algebras but most of them will be too wild  to
fit recognisably into a geometric picture. In the early 1980's most
papers on noncommutative geometry focussed on the noncommutative
torus as the main noncommutative example. Quantum geometry today
contains $q$-planes, $q$-spheres,
$q$-groups $G_q$ (the coordinate rings of $U_q(g)$), $q$-monopoles,
$q$-Lie algebras, $q$-vector fields etc. 

Axiom~2 involves a slight abuse of notation -- what is
quantisation if not a process starting with a Poisson manifold?
The idea is that we need deeper more intrinsically algebraic notions
for the construction of quantum algebras, from which Poisson
manifolds can sometimes (though not necessarily) be obtained by
`classicalisation' with respect to a choice of generators (cf. Lie
algebra contraction). Quantum geometry at present supplies two; one
is the idea of factorisation and the other is the idea of R-matrix or
quasitriangular structures.

Axiom~3 is a more novel issue, usually overlooked even in a Poisson
geometric setting: in classical geometry we demand commutativity {\em
uniformly} for all coordinate algebras. When we relax it for one
geometrical object (e.g. position space), should we not relax it for
another (e.g. momentum space)? And there are many different
`directions' in which one can relax commutativity (e.g. many choices
of Poisson structure or R-matrix) for each object and we need to
choose these consistently. Our quantum spheres have to be consistent
with our quantum planes etc. Section~4 explains how this can be
achieved by means of braid statistics (the braided approach to
q-deformation).

Axiom~4 probably needs the most explanation. A modern way to
think about wave-particle duality is in terms of Fourier theory. The
dual group to position space $\R$ (the group of irreducible
representations) is again $\R$, the momentum space. On the one hand,
points $x$ are fundamental `particles' while on the other hand waves
or points $p$ in momentum space are fundamental `particles'; the two
points of view are related by Fourier transform. When position space
is curved, e.g.
a simple Lie group $G$, then the irreducible representations $\hat G$
do not form a group. However, non-Abelian Fourier theory is still
possible with the right mathematical generalisation. Note that $G$
itself is `geometrical' while the `points' of $\hat G$ are more like
quantum states, so this generalised Fourier theory is an example of a
{\em quantum-geometry transformation}. We would like a similar
elucidation for some class of `group' objects and their duals in any
quantum geometry. Ideally, we might hope for a stronger form of
Axiom~4, which we call the the {\em principle of
representation-theoretic self-duality}: quantum geometry should be
general enough that the duals of `group' objects are again `group'
objects.  This is the case for the kind of quantum geometry coming of
quantum groups and braided groups in Sections~3,4. Moreover, some
groups, like $\R^n$, are self-dual. The self-dual `groups' in quantum
geometry likewise occupy a special place as the simplest geometrical
objects.

Actually, we have  argued\cite{Ma:pri} that the bifurcation into
`geometrical' and `quantum' ideas (dual to each other and related by
a generalised Fourier or quantum geometry transformation) has its
origin in our concept of physical reality itself (the dual nature of
measurement and object being measured), and is not really tied to
groups. More general manifolds with `geometrical' structure have more
complex dual notions of `representation', etc. and self-duality will
pick out geometries occupying a special place. We have postulated
this self-duality
constraint for the quantum geometry of phase space as the
philosophical origin of something like Einstein's
equation\cite{Ma:pri}. 

\section{Elements of quantum geometry}

The present approach to quantum geometry is motivated particularly
from the duality Axiom~4, which leads one to formulate group objects
as Hopf algebras. A Hopf algebra is a unital algebra $H$ and algebra
homomorphisms $\Delta:H\to H\tens H$ (coproduct), $\eps:H\to\C$
(counit) forming a counital coalgebra. The axioms of a counital
coalgebra are just those of a unital algebra with all arrows reversed
(think of the unit element as a map $\C\to H$). So $H^*$ is also a
unital algebra by dualising $\Delta,\eps$. There is also a kind of
`linearised inversion' $S:H\to H$ called the antipode. 

In a suitable setting, if $H$ is a Hopf algebra then essentially
$H^*$ is another Hopf algebra by dualisation. The product of one
determines the coproduct of the other and vice-versa. So Axiom~4 is
satisfied in the strong form. Moreover, if $G$ is a group then its
coordinate ring $\C(G)$ is a Hopf algebra, at least in nice cases. In
the finite case we mean all functions $f$ on $G$ with pointwise
product and $\Delta f=f(\ \cdot\ )$, where $\cdot$ is the group
product. The blanks on the right hand side indicate a function of two
variables, i.e. an element of $\C(G)\tens\C(G)$. So Hopf algebras
generalise the notion of usual groups. Dually paired to $\C(G)$ is
the group algebra $\C G$ (the linear extension of $G$ with $\Delta
g=g\tens g,\eps g=1$) in the finite case or enveloping algebra $U(g)$
in the Lie group case with Lie
algebra $g$. Fourier theory is a linear isomorphism $H\to H^*$ or
$\C(G)\to
\C G$
in the finite case. When $G$ is Abelian, we have $\C G=\C(\hat G)$,
the coordinate ring of the dual group. But when $G$ is non-Abelian, 
$\C G$ is not commutative, so then $\hat G$ only exists as a quantum
geometry with, by definition, coordinate ring $\C G$. One should
consider any noncommutative Hopf algebra as, by
definition, a quantum group.  

\subsection{Toy models of quantum-gravity}

Here we will see how the structural considerations   in Section~2
can force one to concrete Planck scale dynamics. For our
discussion, all quantum geometries that we consider will be `group'
objects, i.e. Hopf algebras, but   the ideas could ultimately be
applied more generally.

Suppose that we fix the position Hopf algebra $H_1$ and momentum
Hopf algebra $H_0$. Instead of Poisson brackets or other guesswork
about the quantum phase space (i.e. instead of guessing
position-momentum
commutation relations) let us proceed structurally and consider all
possible extensions
\eqn{ext}{H_1\to E\to H_0.} 
This means a Hopf algebra $E$ and Hopf algebra maps as shown, obeying
certain conditions\cite{Ma:book}. Then theory tells us that $E$ will
be a cocycle bicrossproduct $E=H_1\bicross
H_0$. The simplest case is with trivial cocycles, in some
cohomological sense `close' to the tensor product $H_1\tens H_0$. 
Then the theorem is that $E$ is a cross product by an action of $H_0$
on $H_1$ and a cross coproduct by a coaction of $H_1$ on $H_0$.

For example, let us consider {\em all possible} extensions
\eqn{extxp}{ \C[x]\to E\to \C[p]}
of 1-dimensional position and momentum spaces $\C[x]$ and $\C[p]$.
These are classical groups and hence Hopf algebras, with $\Delta
x=x\tens 1+1\tens x$ and $\Delta p=p\tens 1+1\tens p$. We consider
all possible $E$, i.e. do not build any relations in by hand, and we
consider the Hamiltonian to be fixed in the form $p^2/2m$. Note that
this approach is
more intrinsic (in the spirit of Einstein's equivalence principle)
than keeping the commutation relations in some canonical form but
varying the Hamiltonian.
  
\begin{propos}\cite{Ma:pla} The possible cocycle-free extensions
(\ref{extxp}) are described by two parameters $\hbar,\grav$ and take
the form $E_{\hbar,\grav}=\C[x]\bicross \C[p]$ with cross relations
and coproduct
\[ [x,p]=\imath\hbar(1-e^{-{x\over M\grav}}),\quad \Delta p=p\tens
e^{-{x\over M\grav}}+1\tens p,\]
where $M$ is a convenient fixed constant of mass dimension.
\end{propos}

The free-fall Hamiltonian $p^2/2m$ gives $\dot
x=v_\infty(1-e^{-{x\over M\grav}})$ to lowest order
in $\hbar$, which can be compared with infalling radial coordinates
$\dot x=-(1-(1+{x\over
2M\grav})^{-1})$ near a black hole of mass $M$. So the parameter
$\grav$ plays a role in
this simple model similar to the gravitational coupling constant.
Also, the particle moves more
and more slowly as it approaches the origin and takes an infinite
time to reach it. Yet when $\grav$ is small, the commutation
relations are $[x,p]=\imath\hbar$ at least for states where one can
say that $x>0$, i.e. away from the origin. Further analysis of this
model gives the effective scales for these gravitational and quantum
limits as $mM>> m_P^2$ and $mM<< m_P^2$, where
$m_P=\sqrt{\hbar/\grav}$ is the Planck mass
\cite{Ma:pla}\cite{Ma:book}.

The same ideas work when the position and momenta are curved. Let
$G_1$ be a Lie group and $g_0$ a Lie algebra with group
$G_0$. One can consider
cocycle-free extensions
\eqn{extG}{ \C(G_1)\to E\to U(g_0).}
The possible extensions turn out\cite{Ma:phy} to correspond
 essentially to solutions of the factorisation problem: Lie
groups factorising into $G_0G_1$. For example, the complexification
of any compact real form $G_0$ of a simple Lie group factorises as
$G_0G_1$ for a certain solvable $G_1$, and gives a
corresponding $E$. The natural Hamiltonian is the quadratic
the Casimir of $g_0$ and induces quantum dynamics on $G_1$ as
position space. 

For example, the quantum group $E=\C(\R^2\lcross \R)\bicross U(so_3)$
corresponds to the Iwasawa factorisation of $SL(2,\C)$. One can
insert two free parameters $\hbar,\grav$ as well. Then $E$ is
generated by the coordinates $x_i$ of $\R^2\lcross\R$ and $e_i$ of
$su_2$, with\cite{Ma:book}
\ceqn{bicso3}{[e_i,e_j]=\imath\hbar \eps_{ijk}e_k,\quad
[e_i,x_j]=\imath\hbar\eps_{ijk}x_k-{\imath\hbar\over
2M\grav}\eps_{ij3}x\cdot
x (1+ {x_3\over M\grav})^{-1}\\
{}[x_i,x_j]=0,\quad\Delta x_i=x_i\tens 1+ (1+{x_3\over M\grav})\tens
x_i\\ 
\Delta e_i=e_i\tens (1+ {x_3\over M\grav})^{-1}+{1\over
M\grav}e_3\tens
x_i(1+{x_3\over M \grav})^{-1}+1\tens e_i.}

By thinking of the $x_i$ as momenta $p_i$ (wave-particle duality
again), one can also consider this $E$ as some kind of
deformation of
$U(\R^3\lcross so_3)$, i.e. of the Poincar\'e enveloping algebra in 3
dimensions. So on the one hand, $E$ is a quantisation of particles on
orbits in $\R^3$ exhibiting singular dynamics, and on the other it is
a deformation of a symmetry algebra. The Minkowski spacetime version
of this model  is very similar and of independent interest
\cite{MaRue:bic}.

Moreover, $E^*$ is another quantum phase space. It solves the
extension problem 
\eqn{dext}{H_0^*\to E^*\to H_1^*.}
When position   space is flat as in (\ref{extxp}) then essentially
$\C[x]^*=\C[p]$ (wave particle duality) and $E^*$ also solves  
(\ref{extxp}). Hence it is of the same form as in Proposition~3.1,
i.e. these $E$ are self-dual quantum groups. 
In the curved position space case, the dual quantum group is
$E^*=U(g_1)\cobicross \C(G_0)$. For example, the dual Hopf algebra to
(\ref{bicso3}) describes quantum particles in
$SU_2=S^3$ moving on orbits under $\R^2\lcross \R$. The explicit
orbits and flows in these models are obtained by solving nonlinear
equations. Also, one can
classicalise and obtain Poisson manifolds of which these quantum
groups are quantisation, although they would not be determined
uniquely as such; see \cite{Ma:book}.

This demonstrates our algebraic approach to Planck
scale physics. It is one of the historical origins of noncommutative
(and noncocommutative) Hopf algebras or quantum groups. Recent work
on bicrossproducts is in \cite{BegMa:qua}. 

\subsection{Quasitriangular structures}

It would be remiss not to mention the more famous Drinfeld-Jimbo
quantum groups $U_q(g)$ and their duals $G_q$. They have little, so
far, to do with Planck scale physics (as far as I know),
arising independently in quite a different physical context. They
can, however, be classicalised and hence viewed  (if we want) as
quantisations of a certain Drinfeld-Sklyanin Poisson bracket on the
Lie group of $g$. At this level, there are connections with the
factorisation problem above\cite{Ma:phy}. Also, they again
demonstrate our Axiom~2 that
quantum geometry has its own intrinsic structure. The intrinsic
structure of $U_q(g)$ is that of a {\em quasitriangular} Hopf
algebra\cite{Dri}. It is a Hopf algebra $H$ equipped with a so-called
(by physicists) universal R-matrix $\CR\in H\tens H$. Its image in
any matrix representation obeys the Yang-Baxter or braid relations.
Such generalised symmetry algebras are relevant to the next section.

The intrinsic structure of $G_q$ is therefore that of a {\em dual
quasitriangular} Hopf algebra. This is a Hopf algebra $H$ equipped
with a skew bicharacter $\CR:H\tens H\to \C$ obeying 
\eqn{dqua}{ g\o h\o \CR(h\t,g\t)=\CR(h\o,g\o)h\t g\t }
for all $h,g\in H$, where $\Delta h=h\o\tens h\t$ is our notation for
the coproduct with output in $H\tens H$ (summation omitted). When
$\CR$, $G_q$ are expanded in $\hbar$ with $q=e^{\hbar/2}$, one
obtains a Poisson bracket. But the quantum world is richer. For
example, there are discrete quantum groups possessing such $\CR$.
Thus the axioms for $H,\CR$ carve out a class of quantum groups
defined intrinsically and `close' to being commutative in the sense
(\ref{dqua}) rather than in the conventional sense of quantisation of
a Poisson bracket.  

\section{Elements of braided geometry}

In this section we explain another approach to quantum geometry,
which has so far been applied mostly in flat space (rather than
having direct contact with the Planck scale), but which has the merit
of solving the uniformity Axiom~3. Ultimately, we would like to see
it combined with the ideas in the preceding section. This {\em
braided geometry} involves a new  kind of mathematics in which
information `flows' along braids and tangles much as it flows along
the wiring in
a computer, except that under- and over-crossings of wires are now
nontrivial {\em braiding operators} $\Psi$. In usual mathematics and
computer science one wires outputs of operations into inputs of other
operations without caring about such crossings, i.e. usual
mathematics
is two-dimensional. By contrast, braided calculations, braided
Feynman diagrams etc. truly exist in a three-dimensional space where
calculations take place. Mathematically, we make use of the theory of
braided categories \cite{JoyStr:bra}. The introduction of algebras,
group theory and geometry in   braided categories is due to the
author, e.g.\cite{Ma:exa}\cite{Ma:introp}.

The idea is that in quantum physics there is another kind of
noncommutativity, namely anticommutativity due
to fermionic statistics. This is a noncommutativity of $\tens$
itself. Thus, when independent fermionic systems  must be
exchanged during a manipulation, one uses supertransposition
$\Psi(b\tens c)=(-1)^{|b||c|}c\tens b$, where $|\ |$ is the degree
$0,1$. For example, the supertensor product $B\tens C$ of two
superalgebras involves $\Psi$, with the result that $cb\equiv (1\tens
c)(b\tens 1)=\Psi(c\tens b)=(-1)^{|c||b|}bc$ in $B\tens C$. The idea
of braided geometry is that $\Psi$, and hence the cross relations of
$B\tens C$, can be much more general than this simple $\pm1$ form.
When $\Psi^2$ is not always the identity, one says that the system
has {\em braid statistics}. Thus,

\begin{itemize} 
\item quantum geometry: $\tens$ usual commutative (bosonic) one,
coordinate algebras noncommutative.
\item braided geometry: $\tens$ non commutative (braid statistics),
coordinate algebras may as well be `commutative' in suitably modified
sense.
\end{itemize}

Just as quantum mechanics was created with the realisation that many
construction do not require commutativity of coordinates, braided
geometry is created by a second and equally deep realisation: {\em
many constructions do not require commutativity of the notion of
independence.} In particular, we can take in place of $-1$ a
dimensionless parameter
$q$ or, more generally, an operator $\Psi$ depending on one or more
parameters $q$. Moreover, {\em classical braided geometry 
$\Rightarrow$ quantum geometry} in a  $q$-deformed  sense because
`braided commutative' generally means noncommutative with respect to
the usual $\tens$. Moreover, specifying the braid statistics
specifies such things coherently between every object and every other
object. Quantum groups still play a role:

\begin{propos}cf.\cite{Dri} All objects $B,C$ on which a quantum
group like
$U_q(g)$ (a quasitriangular Hopf algebra) acts acquire braid
statistics $\Psi(b\tens c)=\CR_i.c\tens \CR^i.b$, where
$\CR=\CR^i\tens\CR_i$ is the universal R-matrix or quasitriangular
structure.
\end{propos}

An example is the quantum-braided plane $\C_q^2$ generated by a
vector of coordinates $\vecx=(x,y)$ obeying $yx=qxy$. It has braiding
and braided-coproduct\cite{Ma:poi}
\ceqn{qplane}{\Psi(x\tens x)=q^2x\tens x,\quad \Psi(y\tens
y)=q^2y\tens y,\quad \Psi(x\tens y)= q y\tens x\\
\Psi(y\tens x)=qx\tens y+(q^2-1)y\tens x,\quad \Delta\vecx=\vecx\tens
1+1\tens\vecx.}
There are braided-plane structures for q-Euclidean and q-Minkowski
spaces. They are isomorphic to their duals  
($q$-wave-particle duality). There are also braided matrices $B(R)$
generated by $\vecu=(u^i{}_j)$ with relations $R_{21}\vecu_1
R\vecu_2=\vecu_2
R_{21}\vecu_1 R$ and\cite{Ma:exa}
\eqn{BR}{\Delta  \vecu=\vecu\tens\vecu,\quad  \Psi(\vecu_1\tens
R\vecu_2)=R\vecu_2 R^{-1}\tens\vecu_1 R}
in a compact notation, for any biinvertible matrix $R$ obeying the
Yang-Baxter equations. Their quotients by $q$-determinant and other
relations give braided versions $BG_q$ of the coordinate rings of
simple Lie groups. They are dual to braided versions $BU_q(g)$ of the
enveloping algebras. Here is a remarkable selfduality phenomenon:

\begin{propos}\cite{Ma:introp}\cite{Ma:lie} When $q\ne 1$ one has
essentially $BG_q\isom BU_q(g)$.
\end{propos}

So in the $q\ne 1$ world there is essentially only one $q$-deformed
object for each simple Lie algebra $g$, which has two limits as $q\to
1$. On the left hand side it becomes $\C(G)$ the commutative
coordinate ring. On the right hand side it becomes the enveloping
algebra $U(g)$. Thus these two features of classical mathematics,
conceptually dual to each other, are different scaling limits of one
object.  In a similar way, one finds that $q$-Minkowski space as a
$2\times 2$ braided matrix is isomorphic to the braided enveloping
algebra of a braided Lie algebra version of $su_2\oplus
u(1)$ \cite{Ma:lie}. This is wave-particle duality in a strong form, 
and is only possible when $q\ne 1$. 

The $q$-Poincar\'e and $q$-conformal groups are also obtained
from braided geometry.  With $q$-Minkowski space as an additive
braided group (like the quantum-braided plane above) one has a
braided adjoint action of the braided-coordinates on themselves. This
is not possible when $q=1$ since the adjoint action is then trivial.
However, when $q\ne 1$ it generates the action of $q$-special
conformal transformations\cite{Ma:geo}. The remnant of this as $q\to
1$ is    
\eqn{conf}{ I\circ{\del\over\del x_i}\circ I=\lim_{q\to
1}{\Ad_{x_i}\over
q-q^{-1}}}
where $I$ is conformal inversion. This is a completely new
group-theoretical picture of conformal transformations as adjoint
action, only possible when $q\ne 1$.  

The above approach to $q$-deformation has been developed over
50--60 papers by the author and collaborators since 1989.  It
provides the correct meaning of $q$ as `braid statistics' (rather
than directly related to $\hbar$) and a systematic solution to the
problem of $q$-deforming everything. Moreover, we see that our
familiar $q=1$ world is merely a special limit of a deeper and more
natural $q\ne 1$ geometry.

\section*{Acknowledgments} I would like to thank the organizers,
H.-D. Doebner, P. Kramer, W. Scherer, V. Dobrev and others for an
outstanding conference.


\end{document}